\documentclass[12pt]{amsart}
\usepackage{amssymb}

\oddsidemargin \evensidemargin
   \title[Combinatorial Formula]
               {A combinatorial formula for non-symmetric
                         Macdonald polynomials}

   \author{J. Haglund}
   \thanks{Work supported by NSA grant MSPF-02G-193 (J.H.)}

   \author{M. Haiman}
   \thanks{Work supported by NSF grant DMS-0301072 (M.H.)}

   \author{N. Loehr}
   \thanks{Work supported by NSF postdoctoral research fellowship 
               DMS-0303178 (N.L.)}

   \address[J.H.]{Dept.\ of Mathematics\\
            University of Pennsylvania\\
            Philadelphia, PA}

   \address[M.H.]{Dept.\ of Mathematics\\
            University of California\\
            Berkeley, CA}

   \address[N.L.]{Dept.\ of Mathematics\\
            The College of William and Mary\\
            Williamsburg, VA}

   \email[J.H.]{jhaglund@math.upenn.edu}
   \email[M.H.]{mhaiman@math.berkeley.edu}
   \email[N.L.]{nick@math.wm.edu}

   \date{February 9, 2007}

\newtheorem{thm}{Theorem}[subsection]
\newtheorem{lemma}[thm]{Lemma}
\newtheorem{prop}[thm]{Proposition}
\newtheorem{cor}[thm]{Corollary}

\theoremstyle{definition}
\newtheorem{defn}[thm]{Definition}

\theoremstyle{remark}

\newtheorem*{remarks}{Remarks}

\DeclareMathOperator{\diagram}{dg}
\newcommand{\skydg}{\diagram'}
\newcommand{\augdg}{\widehat{\diagram}}
\DeclareMathOperator{\arm}{\text{arm}}
\DeclareMathOperator{\leg}{\text{leg}}
\DeclareMathOperator{\Des}{Des}
\DeclareMathOperator{\Inv}{Inv}

\DeclareMathOperator{\maj}{maj}

\DeclareMathOperator{\inv}{inv}
\DeclareMathOperator{\coinv}{coinv}

\newcommand{\Acal}{{\mathcal A}}
\newcommand{\Ecal}{{\mathcal E}}
\newcommand{\Hcal}{{\mathcal H}}

\newcommand{\NN}{{\mathbb N}}
\newcommand{\QQ}{{\mathbb Q}}
\newcommand{\ZZ}{{\mathbb Z}}

\newcommand{\defeq}{\underset{\text{{\it def}}}{=}}

\newlength{\cellsize}
\newcommand\tableau[1]{
\vcenter{
\let\\=\cr
\baselineskip=-16000pt
\lineskiplimit=16000pt
\lineskip=0pt
\halign{&\tableaucell{##}\cr#1\crcr}}}

\newcommand{\tableaucell}[1]{{%
\def \arg{#1}\def \void{}%
\ifx \void \arg
\vbox to \cellsize{\vfil \hrule width \cellsize height 0pt}%
\else
\unitlength=\cellsize
\begin{picture}(1,1)
\put(0,0){\makebox(1,1){$#1$}}
\put(0,0){\line(1,0){1}}
\put(0,1){\line(1,0){1}}
\put(0,0){\line(0,1){1}}
\put(1,0){\line(0,1){1}}
\end{picture}%
\fi}}

\begin{document}


\begin{abstract}
We give a combinatorial formula for the non-symmetric Macdonald
polynomials $E_{\mu}(x;q,t)$.  The formula generalizes our previous
combinatorial interpretation of the integral form symmetric Macdonald
polynomials $J_{\mu}(x;q,t)$.  We prove the new formula by verifying
that it satisfies a recurrence, due to Knop and Sahi, that
characterizes the non-symmetric Macdonald polynomials.
\end{abstract} 

\maketitle

\section{Introduction}
\label{sec:intro}

In a previous paper \cite{HaHaLo05}, we gave a
combinatorial formula for Macdonald's symmetric functions (for root
systems of type $A$).  In fact, we gave two formulas, one for the
transformed Macdonald polynomials $\tilde{H}_{\mu }(x;q,t)$, which
appear in the positivity theorem of \cite{Hai01} and are connected to
the geometry of Hilbert schemes, and another \cite[Proposition
8.1]{HaHaLo05} for the {\it integral forms} $J_{\mu }(x;q,t)$
\cite{Mac95}, which are scalar multiples of the classical monic forms
$P_{\mu }(x;q,t)$, and thus closer to the original viewpoint on
Macdonald polynomials as orthogonal polynomials associated with a root
system \cite{Mac00}.

The symmetric Macdonald polynomials $P_{\mu }(x;q,t)$ can be
reconstructed from a broader theory of {\it non-symmetric Macdonald
polynomials} $E_{\mu }(x;q,t)$, defined and developed by Opdam and
Heckman \cite{Opdam95}, Macdonald \cite{Mac96}, and Cherednik
\cite{Chered95}.  Many important aspects of the symmetric theory are
best understood with the help of extra structure available in the
non-symmetric theory.  These include the norm and evaluation formulas
\cite{Chered95a}, \cite{Chered95b}, as well as the earlier
combinatorial formula of Knop and Sahi \cite{KnSa97} for the {\it Jack
polynomials} $P_{\mu }(x;\alpha ) = \lim_{t\rightarrow 1}P_{\mu
}(x;t^{\alpha },t)$.  In fact, Knop and Sahi's formula was derived
from its non-symmetric counterpart.  It is similarly possible to
derive our formula for $P_{\mu }(x;q,t)$ in \cite{HaHaLo05} from the
formula for $E_{\mu }(x;q,t)$ obtained here---see \S
\ref{subsec:stable}.  Apart from its value as a path to results in the
symmetric theory, the non-symmetric theory has lately taken on a life
of its own, with works such as \cite{Ion04}, \cite{Ion05},
\cite{Mason06} suggesting that the $E_{\mu }(x;q,t)$ for any root
system should possess intrinsic Lie-theoretic significance.  It also
seems to us that the best hope for progress on the problem of
extending the combinatorial theory to other root systems lies in the
non-symmetric setting.

The proof of our formula uses a recurrence discovered independently by
Knop \cite{Kno97} and Sahi \cite{Sah96} and used by them to prove
Macdonald's integrality conjecture.  In retrospect, The Knop--Sahi
recurrence is the special case in type $A$ of Cherednik's intertwiner
formula \cite{Chered97}.

In the next section we fix notation, review the definition of
non-symmetric Macdonald polynomials, and explain the Knop--Sahi recurrence.
In \S \ref{sec:formula} we state our main theorem, the combinatorial
formula.  The proof is in \S \ref{sec:proof}.  In \S
\ref{sec:comparison} we compare the new formula to our earlier one for
the symmetric case.  Finally, we give a table of type $A$
non-symmetric Macdonald polynomials in an appendix, for the
convenience of the reader who might wish to compare examples
calculated from our formula against known values.

\section{Notation and definitions}
\label{sec:defs}

\subsection{Non-symmetric Macdonald polynomials}
\label{subsec:E-mu}

Our notational conventions mostly conform to those in Cherednik's
paper \cite{Chered95} and the expositions by Cherednik \cite[Chapter
3]{Chered05}, Macdonald \cite{Mac03}, and the second author
\cite{Hai06}, particularized to $GL_{n}$ as indicated below.  In
particular, $t$ is the Hecke algebra parameter, and $q$ is the formal
exponential of the null root.

The integer $n$ is fixed throughout.  The weight lattice of $GL_{n}$
is $X = \ZZ ^{n}$; the simple roots are $\alpha _{i} = e_{i}-e_{i+1}$,
where $e_{i}$ is the $i$-th unit vector.  Identifying the co-weight
lattice $X^{\vee }$ with $X$ using the standard inner product on $\ZZ
^{n}$, the simple coroots $\alpha _{i}^{\vee }$ coincide with the
simple roots $\alpha _{i}$.  Hence the dominant weights, {\it i.e.},
$\lambda \in X$ such that $\langle \lambda , \alpha
_{i}^{\vee}\rangle\geq 0$ for all $i$, are the weakly decreasing
sequences $(\lambda _{1}\geq \cdots \geq \lambda _{n})$.  We write
$x^{\lambda } = x_{1}^{\lambda _{1}}\cdots x_{n}^{\lambda _{n}}$ for
the formal exponential of a weight $\lambda $.

The affine weight lattice is $\widehat{X} = X\oplus \ZZ \delta $,
where $\delta $ is the smallest positive imaginary root, or {\it null
root.}  The simple affine roots are $\alpha _{1},\ldots,\alpha _{n}$
as above, and $\alpha _{0} = \delta -\theta $, where $\theta
=e_{1}-e_{n}$ is the highest root of $GL_{n}$.  The set of positive
affine (real) roots is $\widehat{R}_{+} = \{e_{i}-e_{j}+k\delta :
\text{$i\not =j$, $k\in \NN $, and $k > 0$ if $i > j$} \}$.  The
formal exponential $x^\delta $ is denoted by $q$, so $x^{\alpha _{i}}
= x_{i}/x_{i+1}$ for $i\not =0$, and $x^{\alpha _{0}} = qx_{n}/x_{1}$.
The group algebra $\QQ (t) \widehat{X}$ is thereby identified with the
Laurent polynomial ring $\QQ (t)[q^{\pm 1},x_{1}^{\pm
1},\ldots,x_{n}^{\pm 1}]$.  Extending scalars to $\QQ (q,t)$, we
further identify $\QQ (t) \widehat{X}$ with a subring of $\QQ (q,t) X
= \QQ (q,t)[x_{1}^{\pm 1},\ldots,x_{n}^{\pm 1}]$.

{\it Cherednik's inner product} on $\QQ (q,t)X$ is defined by
\begin{equation}
\langle f,g \rangle _{q,t}= [x^{0}](f\overline{g}\Delta _{1}),
\end{equation}
where $\overline{\, \cdot \, }$ is the involution $\overline{q} =
q^{-1}$, $\overline{t} = t^{-1}$, $\overline{x_{i}} = x_{i}^{-1}$, and
$\Delta _{1} = \Delta /([x^{0}]\Delta )$, with
\begin{equation}
\Delta =\prod _{\alpha \in \widehat{R}_{+}}\frac{1-x^{\alpha
}}{1-tx^{\alpha }} = \prod _{i<j} \prod _{k=0}^{\infty }
\frac{(1-q^{k}x_{i}/x_{j})(1-q^{k+1}x_{j}/x_{i})}{(1-tq^{k}x_{i}/x_{j})(1-tq^{k+1}x_{j}/x_{i})}.
\end{equation}
Here $\Delta \in \QQ [[q,t]]$ is a formal power series in $q,t$.  For
such a series $f$, $[x^{\lambda }]f\in \QQ [[q,t]]$ denotes its
coefficient of $x^{\lambda }$, taken term by term.  It is known that
$[x^{\lambda }]\Delta _{1}$ is a rational function of $q,t$ for all
$\lambda $, and that $\overline{\Delta _{1}} = \Delta _{1}$.  Hence
$\langle f,g \rangle _{q,t}\in \QQ (q,t)$, and $\langle g,f \rangle _{q,t}=
\overline{\langle f,g \rangle }_{q,t}$.

The {\it Bruhat order} is the partial ordering $<$ on $X$ induced by
identifying $X$ with the set of minimal coset representatives in
$\widehat{W}/W_{0}$, where $W_{0} = S_{n}$ is the Weyl group of
$GL_{n}$, and $\widehat{W} = W_{0}\ltimes X$ is the extended affine
Weyl group, equipped with its usual Bruhat order.  To make it explicit
for $GL_{n}$, if $i<j$, $\lambda _{i}<\lambda _{j}$ and $\sigma _{ij}$
is the transposition $(i\;j)$, then $\lambda >\sigma _{ij}(\lambda )$,
and if in addition $\lambda _{j}-\lambda _{i}>1$, then $\sigma
_{ij}(\lambda )>\lambda +e_{i}-e_{j}$.  The Bruhat order is the
transitive closure of these relations.  Note that distinct cosets in
$X$ of the root lattice $Q = \ZZ \{\alpha _{1},\ldots,\alpha _{n} \}$
are incomparable in the Bruhat order.

\begin{defn}\label{def:E-mu}
The {\it non-symmetric Macdonald polynomials} $E_{\mu }(x;q,t)\in \QQ
(q,t)X$ $(\mu \in X)$ are uniquely characterized by the conditions
\begin{itemize}
\item [(i)] triangularity: $E_{\mu } \in x^{\mu } + \QQ (q,t)\{x^{\lambda
}:\lambda <\mu \}$,
\item [(ii)] orthogonality: $\langle E_{\lambda },E_{\mu } \rangle _{q,t}=
0$ for $\lambda \not =\mu $.
\end{itemize}
\end{defn}

We remark that Opdam \cite{Opdam95} and Macdonald \cite{Mac96},
\cite{Mac03} use a stronger ordering that is superficially easier to
describe but not as natural, since, as Sahi \cite{Sahi00} observed,
the coefficient of $x^{\lambda }$ in $E_{\mu }(x;q,t)$ in non-zero if
and only if $\lambda \leq \mu $ in Bruhat order.

The proof that the polynomials $E_{\mu }$ exist uses an affine Hecke
algebra representation which is also needed in order to state the
Knop--Sahi recurrence, and which we now recall.  The {\it affine Hecke
algebra} is the $\QQ (t)$-algebra $\Hcal $ with generators $T_{i}$
($i=0,\ldots,n-1$) satisfying the braid relations
\begin{equation}
\begin{aligned}
T_{i}T_{i+1}T_{i} & = T_{i+1}T_{i}T_{i+1}\\
T_{i}T_{j} & = T_{j}T_{i}\qquad (i-j\not =\pm 1),
\end{aligned}
\end{equation}
where all indices are modulo $n$, and the quadratic relations
\begin{equation}\label{e:quadratic}
(T_{i}-t)(T_{i}+1)=0.
\end{equation}
Here we depart slightly from \cite{Chered05}, \cite{Hai06},
\cite{Mac03}, which take $(T_{i}-t^{1/2})(T_{i}+t^{-1/2})=0$ instead.
We will not use the parameter $t^{1/2}$.

The (unextended) affine Weyl group $W_{a}$ is a Coxeter group with
generators $s_{i}$ ($i=0,\ldots,n-1$) satisfying the above braid
relations.  They act on $\widehat{X}$, and by extension on $\QQ
(t)\widehat{X}$ and $\QQ (q,t)X$, as affine reflections
\begin{equation}
s_{i}(\lambda
)=\lambda -\langle \lambda ,\alpha _{i}^{\vee } \rangle \alpha _{i}.
\end{equation}
For $i\not =0$, $s_{i}$ is just the transposition $\sigma _{i,i+1}$,
while $s_{0}$ is given explicitly by
\begin{equation}
s_{0}f(x_{1},\ldots,x_{n}) = f(qx_{n},x_{2},\ldots,x_{n-1},x_{1}/q).
\end{equation}
{\it Cherednik's representation} of $\Hcal $ on $\QQ (q,t)X$ is
defined by the formula
\begin{equation}\label{e:cherednik}
T_{i}\, x^{\lambda } = t\, x^{s_{i}(\lambda )}+(t-1)\frac{x^{\lambda
}-x^{s_{i}(\lambda )}}{1-x^{\alpha _{i}}},\quad \text{for all
$i=0,\ldots,n-1$},
\end{equation} 
where for $i=0$, we recall that $x^{\alpha _{0}} = qx_{n}/x_{1}$.

Denote the standard basis elements of $\Hcal $ by $T_{w}$ ($w\in
W_{a}$).  For $\tau (\beta )\in W_{a}$ corresponding to translation by
a dominant weight $\beta \in Q^{\vee } = Q$, set $Y^{\beta
}=t^{-\langle \beta, \rho ^{\vee } \rangle}T_{\tau (\beta )}$, acting
on $\QQ (q,t)X$ via Cherednik's representation.  Here $\rho ^{\vee } =
(\sum _{\alpha \in R^{+}}\alpha ^{\vee })/2$, as usual.  The operators
$Y^{\beta }$ commute, are unitary with respect to Cherednik's inner
product, and are lower triangular with respect to Bruhat order on the
basis $\{x^{\lambda }:\lambda \in X \}$.  The polynomials $E_{\mu
}(x;q,t)$ are their simultaneous eigenfunctions.  Bernstein's
relations \cite[(4.2.4)]{Mac03}
%
%
in $\Hcal $ and a computation of the eigenvalues of $Y^{\beta }$ imply
the relations
\begin{equation}\label{e:intertwiners}
E_{s_{i}(\mu )}(x;q,t) = \left(T_{i} + \frac{1-t}{1-q^{\langle \mu ,\alpha
_{i}^{\vee } \rangle} t^{\langle w_{\mu }(\rho ),\alpha _{i}^{\vee
} \rangle}} \right)E_{\mu }(x;q,t),\quad (i\not =0)
\end{equation}
for $\mu _{i}>\mu _{i+1}$ (that is, $s_{i}(\mu )>\mu $), where $w_{\mu
}\in W_{0}$ is the {\it maximum-length} permutation such that $w_{\mu
}^{-1}(\mu )$ is dominant, and $\rho =(n-1,n-2,\ldots,0)$.  A more detailed
argument can be found in \cite[Theorem 4.2b]{Kno97}.

Cherednik's intertwiner theory yields a version of
\eqref{e:intertwiners} for $i=0$, but for $GL_{n}$ it is simpler to
make use of the symmetry of the affine root system.  Define
automorphisms $\pi $ of $X$ and $\Psi $ of $\QQ (q,t)X$ by
\begin{equation}\label{e:pi-Psi}
\begin{gathered}
\pi (\lambda _{1},\ldots,\lambda _{n}) = (\lambda _{n}+1,\lambda
_{1},\ldots,\lambda _{n-1})\\
\Psi f(x_{1},\ldots,x_{n}) = x_{1}f(x_{2},\ldots,x_{n},q^{-1}x_{1}).
\end{gathered}
\end{equation}
If we define a linear automorphism $\phi $ of $\widehat{X}$ by $\phi
(\lambda _{1},\ldots,\lambda _{n}) = (\lambda _{n},\lambda
_{1},\ldots,\lambda _{n-1})-\lambda _{n}\delta $, $\phi \delta =\delta
$, then $\Psi f = x_{1}\Psi 'f$, where $\Psi '(x^{\lambda }) = x^{\phi
(\lambda )}$.  Now, $\phi $ stabilizes $\widehat{R}_{+}$, hence $\Psi
'$ fixes $\Delta $ and is an isometry of Cherednik's inner product
$\langle \cdot ,\cdot \rangle _{q,t}$, which in turn implies that
$\Psi $ is an isometry of $\langle \cdot ,\cdot \rangle _{q,t}$.
Moreover, $\Psi (x^{\lambda }) = q^{-\lambda _{n}}x^{\pi (\lambda )}$,
and $\pi $ preserves the Bruhat order.  It therefore follows from
Definition~\ref{def:E-mu} that
\begin{equation}\label{e:Psi}
E_{\pi (\mu )}(x;q,t) =  q^{\mu _{n}} \Psi E_{\mu }(x;q,t).
\end{equation}
Equations \eqref{e:intertwiners} and \eqref{e:Psi} are the Knop--Sahi
recurrence.

It is immediate from the definition or from \eqref{e:Psi} that if $\mu
'=\mu +(r,r,\ldots,r)$, then $E_{\mu '} = (x_{1}\cdots
x_{n})^{r}E_{\mu }$.  Therefore, without loss of generality, we may
restrict attention to weights $\mu \in \NN ^{n}$, and we will do so
from now on.  Then $\mu $ is a {\it composition}, and dominant weights
are {\it partitions}.

\begin{lemma}\label{lem:Knop-recurrence}
The Macdonald polynomials $E_{\mu }$ for $\mu \in \NN ^{n}$ are
uniquely determined by the the initial value $E_{{\mathbf 0}}= 1$,
together with
\begin{itemize}
\item [(a)] equation \eqref{e:Psi}, and
\item [(b)] the special case of \eqref{e:intertwiners} in which $\mu
_{i+1}=0$.
\end{itemize}
\end{lemma}

\begin{proof}
If $\mu _{1}>0$, then $\mu =\pi (\nu )$, where again $\nu \in \NN
^{n}$.  By induction on the sum of the parts, we can assume $E_{\nu }$
already determined, and apply (a).  If $\mu _{1}=0$ and $\mu _{j}>0$
for some $j$, we can reduce to the case $\mu _{1}>0$ by repeated
applications of (b).
\end{proof}

We should point out that notational conventions in earlier literature
sometimes differ from ours.  Ion \cite{Ion04}, \cite{Ion05} and Marshall
\cite{Mar99} change $q$, $t$ to $q^{-1}$, $t^{-1}$, and our $t$ is
$q^{k}$ in Macdonald \cite{Mac96}.  Both Knop \cite{Kno97} and
Marshall \cite{Mar99} reverse the indexing of the variables
$x_{1},\ldots,x_{n}$, so for them, monomials $x^{\lambda }$ with
$\lambda _{1}\geq \cdots \geq \lambda _{n}$ correspond to {\it
anti-dominant} weights.  The equation $\Phi\xi_1=q\xi_n\Phi$ in
\cite[\S 4]{Kno97} contains a typographical error, and should instead
read $\Phi\xi_1=q^{-1}\xi_n\Phi$.

\subsection{Diagrams}
\label{subsec:diagram-def}

We visualize a composition $\mu=(\mu_1,\ldots,\mu_n)\in\NN^n$ as a
diagram consisting of $n$ columns, with $\mu_i$ boxes in column $i$.
Formally, the \emph{column diagram} of $\mu$ is the set
\begin{equation}
\skydg(\mu)=\{ (i,j)\in\NN^2: 1\leq i\leq n,\ 1\leq j\leq\mu_i \}
\end{equation}
in Cartesian coordinates, the abscissa $i$ indexing the columns, and
the ordinate $j$ indexing the rows.  The prime serves as a reminder
that the parts of $\mu$ are the columns rather than the rows,
transpose to the usual conventions for partition diagrams.

For example, $\mu=(2,1,3,0,0,2)$ is represented by the diagram
\begin{equation}
\skydg (\mu )= \tableau{  &  &{}&  &  &  \\ 
                        {}&  &{}&  &  &{}\\ 
                        {}&{}&{}& & &{}}\; .
\end{equation}
We will also need the \emph{augmented diagram} of $\mu$, defined by
\begin{equation}
\augdg(\mu)=\skydg(\mu)\cup\{ (i,0): 1\leq i\leq n \},
\end{equation}
{\it i.e.}, we adjoin $n$ extra boxes in row $0$, thus adding a box to
the bottom of each column.

\subsection{Arms and legs}
\label{subsec:arm-leg}

Given $\mu\in\NN^n$ and $u=(i,j)\in\skydg(\mu)$, define
\begin{equation}\label{e:arms}
\begin{aligned}
\leg (u) & = \{(i,j')\in \skydg(\mu) : j'>j \}\\
\arm ^{\text{left}}(u) & = \{(i',j)\in \skydg(\mu): i'<i,\, \mu _{i'}\leq \mu
_{i} \}\\
\arm ^{\text{right}}(u) & = \{(i',j-1)\in \augdg(\mu) : i'>i,\, \mu
_{i'}< \mu _{i} \}\\ 
\arm (u) & = \arm ^{\text{left}}(u) \cup \arm ^{\text{right}} (u),
\end{aligned}
\end{equation}
and set
\begin{equation}\label{e:a(u)-l(u)}
\begin{gathered}
l(u) = | \leg (u) | = \mu _{i}-j\\
a(u) = | \arm (u) |.
\end{gathered}
\end{equation}

For example, for $\mu=(3,1,2,4,3,0,4,2,3)$ and $u=(5,2)$,
the cells belonging to $\leg(u)$, $\arm^{\text{left}}(u)$, and 
$\arm^{\text{right}}(u)$ are marked by $x$, $y$, $z$ in the following figure:
\begin{equation}
\skydg(\mu)=
\tableau{   &   &   &{ }&   &   &{ }&   &   \\
         { }&   &   &{ }&{x}&   &{ }&   &{ }\\
         {y}&   &{y}&{ }&{u}&   &{ }&{ }&{ }\\
         { }&{ }&{ }&{ }&{ }&   &{ }&{z}&{ }
}\; ,
\end{equation}
giving $l(u)=1$, $a(u)=3$.  Note that $\arm (u)$ is a subset of the
augmented diagram, and may contain boxes in row $0$, if $u$ is in row
$1$.  When $\mu $ is anti-dominant, $\skydg (\mu )$ is a reversed
partition diagram, all right arms are empty, and the leg and left arm
reduce to the usual notions for partition diagrams.  We remark that
the preceding definitions agree with those in \cite{Kno97},
\cite{KnSa97}, after a suitable change of coordinates, and exchanging
arms with legs.

With these definitions, if $\mu _{i}>\mu _{i+1}$, and $u=(i,\, \mu
_{i+1}+1)$, then \eqref{e:intertwiners} takes the form
\begin{equation}\label{e:arm-leg-intertwiners}
E_{s_{i}(\mu )}(x;q,t) = \left(T_{i} + \frac{1-t}{1-q^{l(u)+1}
t^{a(u)}} \right)E_{\mu }(x;q,t).
\end{equation}

\subsection{Integral forms}
\label{subsec:integralform-def}
For $\mu\in\NN^n$, the \emph{integral form} non-symmetric 
Macdonald polynomials are defined by
\begin{equation}\label{e:Ecal}
\Ecal _{\mu }(x;q,t) = \prod _{u\in \skydg (\mu )} \left(
1-q^{l(u)+1}t^{a(u)+1} \right) E_{\mu }(x;q,t).
\end{equation}
Knop \cite[Corollary 5.2]{Kno97} proved that $\Ecal_{\mu}$ has
coefficients in $\ZZ [q,t]$, as will also be seen from our formula.

\section{The combinatorial formula}
\label{sec:formula}

In this section we present our main result, Theorem \ref{thm:main}, a
combinatorial formula for the non-symmetric Macdonald polynomials
$E_{\mu }$.  The formula is a sum whose terms are indexed by
combinatorial objects $\sigma $ called {\it fillings of $\mu $}.  Each
term in the sum depends on a pair of combinatorial statistics $\maj
(\widehat{\sigma })$ and $\coinv (\widehat{\sigma })$, which are
variants of their counterparts defined in \cite{HaHaLo05},
\cite{HaHaLo05a} for the symmetric case.  We begin by defining these
combinatorial data.

\subsection{Fillings}
\label{subsec:fillings}

A {\it filling of $\mu $} is a function 
\begin{equation}
\sigma \colon \skydg(\mu)\rightarrow [n],
\end{equation}
where $[n] = \{1,\ldots,n \}$.  The associated {\it augmented filling}
is the filling $\widehat{\sigma }\colon \augdg(\mu)\rightarrow [n]$ of
the augmented diagram such that $\widehat{\sigma }$ agrees with
$\sigma $ on $\skydg(\mu)$, and $\widehat{\sigma }((j,0)) = j$ for $j
= 1,\ldots,n$.  Distinct lattice squares $u,v\in \NN ^{2}$ are
said to {\it attack} each other if either
\begin{itemize}
\item [(a)] they are in the same row, {\it i.e.}, they have the form
$(i,j), (i', j)$, or
\item [(b)] they are in consecutive rows, and the box in the lower row
is to the right of the one in the upper row, {\it i.e.}, they have the
form $(i,j), (i',j-1)$ with $i<i'$.
\end{itemize}
A filling $\widehat{\sigma }\colon \augdg(\mu)\rightarrow [n]$ is {\it
non-attacking} if $\widehat{\sigma }(u)\not = \widehat{\sigma }(v)$
for every pair of attacking boxes $u,v\in\augdg(\mu)$.  We say that a
filling $\sigma \colon \skydg (\mu )\rightarrow [n]$ of $\mu $ is {\it
non-attacking} if its associated augmented filling $\widehat{\sigma }$
is non-attacking (non-attacking fillings are called {\it admissible}
in \cite{KnSa97}).

\subsection{Descents and major index}
\label{subsec:des-maj}

We write $d(u) = (i,j-1)$ for the box directly below a box $u =
(i,j)$.  For any lattice-square diagram $S\subseteq \NN ^{2}$,
a {\it descent} of a filling $\widehat{\sigma}\colon S\rightarrow [n]
$ is a box $u\in S$ such that $d(u)\in S$ and $\widehat{\sigma}
(u)>\widehat{\sigma} (d(u))$.  Note that the descents of a filling
$\widehat{\sigma }$ of $\augdg (\mu )$ are contained in $\skydg (\mu
)$.  For such a filling, we define
\begin{align}
\Des (\widehat{\sigma }) &= \{\text{descents of $\widehat{\sigma }$} \},\\
\maj (\widehat{\sigma }) &= \sum _{u\in \Des (\widehat{\sigma })} (l(u)+1).
\end{align}

\subsection{Inversions}\label{subsec:inv}

The {\it reading order} on a lattice-square diagram $S$ is the total
ordering of the boxes in $S$ row by row, from top to bottom, and from
{\it right to left} within each row.  In symbols, $(i,j)<(i',j')$ if
$j>j'$, or if $j=j'$ and $i>i'$.  An {\it inversion} of a filling
$\widehat{\sigma }\colon S\rightarrow [n]$ is a pair of boxes $u,v\in
S$ such that $u$ and $v$ attack each other, $u < v$ in the reading
order, and $\widehat{\sigma }(u)>\widehat{\sigma }(v)$.  When $\sigma 
$ is a filling of $\mu $ and $\widehat{\sigma }$ is its associated
augmented filling, we define
\begin{align}
\Inv (\widehat{\sigma }) &= \{\text{inversions of $\widehat{\sigma }$} \},\\
\inv (\widehat{\sigma}) &= | \Inv (\widehat{\sigma }) | - |\{i<j:\mu
_{i}\leq \mu _{j}
\}|- \sum _{u\in \Des (\widehat{\sigma })} a(u), \label{e:inv} \\
\coinv (\widehat{\sigma }) & = \bigl(\sum _{u\in \skydg (\mu )} a(u)
\bigr) - \inv (\widehat{\sigma }). \label{e:coinv}
\end{align}
Note that every pair of boxes in row $0$ counts as an inversion of
$\widehat{\sigma }$, although the second term in the expression for
$\inv (\widehat{\sigma })$ has the effect of cancelling part of this
contribution.

\subsection{Example}
\label{subsec:example}

The figure below shows the augmentation $\widehat{\sigma }$ of a
non-attacking filling $\sigma $ of $\mu =(2,1,3,0,0,2)$, and the
arm-lengths $a(u)$ for each box $u\in \skydg (\mu )$.
\begin{equation}\label{e:example}
\widehat{\sigma } = \;\tableau{
&		&	2\\
{\mathbf 6}&	&	{\mathbf 4}&	&	&	5\\
1&		2&	3&		&	&	5\\
1&		2&	3&		4&	5&	6
}\;, \qquad a(u) = \;  \tableau{
&	&	1\\
1&	&	2&	&	&	1\\
3&	2&	5&	&	&	2\\
{}&	{}&	{}&	{}&	{}&	{}
}\;.
\end{equation}
The bottom row is row $0$.  The descent set $\Des (\widehat{\sigma })$
consists of the two boxes marked in boldface, giving $\maj
(\widehat{\sigma }) = 3$.  There are $\binom{6}{2} = 15$ inversions in
row $0$, $\binom{4}{2} = 6$ in row $1$, and $1$ in row $2$, plus $3$
inversions between rows $2$ and $1$, for a total of $25$ inversions,
giving $\inv (\widehat{\sigma }) = 25 - 7-3 = 15$, $\coinv
(\widehat{\sigma }) = 17-\inv (\widehat{\sigma })=2$.

\subsection{The formula}
\label{subsec:main}

We can now state our main result.

\begin{thm}\label{thm:main}
The non-symmetric Macdonald polynomials $E_{\mu }$ are given by the
formula
\begin{equation}\label{e:main}
E_{\mu }(x;q,t) = \sum _{\substack{\sigma \colon \mu \rightarrow [n]\\
\text{non-attacking}}} x^{\sigma } q^{\maj (\widehat{\sigma })} t^{\coinv
(\widehat{\sigma })} \prod _{\substack{u\in \skydg (\mu )\\
\widehat{\sigma }(u)\not =\widehat{\sigma }(d(u))}}
\frac{1-t}{1-q^{l(u)+1}t^{a(u)+1}},
\end{equation}
where $x^{\sigma } = \prod _{u\in \skydg (\mu )} x_{\sigma (u)}$.
\end{thm}

\begin{cor}\label{cor:integral-form}
The integral form non-symmetric Macdonald polynomials $\Ecal _{\mu }$
in \eqref{e:Ecal} are given by
\begin{equation}\label{e:main-integral}
\Ecal _{\mu }(x;q,t) = \sum _{\substack{\sigma \colon \mu \rightarrow [n]\\
\text{non-attacking}}} x^{\sigma } q^{\maj (\widehat{\sigma })} t^{\coinv
(\widehat{\sigma })} \prod _{\substack{u\in \skydg (\mu )\\
\widehat{\sigma }(u) =\widehat{\sigma }(d(u))}}
(1-q^{l(u)+1}t^{a(u)+1}) \prod _{\substack{u\in \skydg (\mu )\\
\widehat{\sigma }(u)\not =\widehat{\sigma }(d(u))}} (1-t).
\end{equation}
\end{cor}

Let us make some remarks about the combinatorial nature of
\eqref{e:main} and \eqref{e:main-integral}.  It follows trivially from
\eqref{e:intertwiners} and \eqref{e:Psi} that there exist expressions
for $E_{\mu }$ as a sum of terms of the form
\begin{equation}
\pm x^{\lambda }q^{c}t^{d}\prod _{(a,b)\in S}\frac{1-t}{1-q^{a}t^{b}}.
\end{equation}
In fact, using Cherednik's intertwiners, the same conclusion can be
drawn for non-symmetric Macdonald polynomials associated with any root
system.  The reader might think at first that in \eqref{e:main} we
have merely organized the terms in such an expression, but this is not
so.

For one thing, \eqref{e:main} and \eqref{e:main-integral} are {\it
positive} formulas in at least two senses.  First, if the parameters
are specialized to real numbers $0<q,t<1$ then \eqref{e:main}
becomes a polynomial in $x$ with positive coefficients.  Second, if we
set $q=t^{k}$ for any integer $k\geq 0$, then \eqref{e:main-integral}
shows that
\begin{equation}
\Ecal _{\mu }(x;t^{k},t)/(1-t)^{|\mu |}
\end{equation}
is a polynomial in $x$ and $t$ with positive coefficients.

A more subtle point is that \eqref{e:main} is a {\it canonical}
formula.  In order to apply \eqref{e:intertwiners} and \eqref{e:Psi}
directly, one must first choose a reduced factorization in the affine
Weyl group for the minimal element $w$ that represents $\mu \in X\cong
\widehat{W}/W_{0}$, and the resulting expression for $E_{\mu }$ will
depend on this choice.  Formula \eqref{e:main} involves no such
auxiliary choices.  This is not to say, however, that \eqref{e:main}
is the unique, or even the simplest, positive combinatorial formula
for $E_{\mu }$.  For example, if $\mu = (r,r,\ldots,r)-(\nu
_{n},\ldots,\nu _{1})$, it is easy to deduce from the definition the
following identity:
\begin{equation}\label{e:E-mu-complement}
E_{\mu }(x_{1},\ldots,x_{n};q,t) = (x_{1}\cdots x_{n})^{r}E_{\nu
}(x_{n}^{-1},\ldots,x_{1}^{-1};q,t).
\end{equation}
It may happen that $\nu $ has fewer non-attacking fillings than $\mu
$, in which case \eqref{e:E-mu-complement} combined with formula
\eqref{e:main} for $E_{\nu }$ yields a simpler formula for $E_{\mu }$
than \eqref{e:main}.  It would be interesting to know whether there is
some other combinatorial formula that would manifest the symmetry in
\eqref{e:E-mu-complement} directly.

\subsection{Inversion triples}
\label{subsec:inv-triples}

As in \cite{HaHaLo05}, we can reformulate $\inv (\widehat{\sigma })$
as the number of suitably defined {\it inversion triples}.  The
definition \eqref{e:inv} and the reformulation below will both be used
in the proof of Theorem~\ref{thm:main} in \S \ref{sec:proof}.  We have
placed the reformulation in this section because it clarifies the
nature of the statistics $\inv (\widehat{\sigma })$ and $\coinv
(\widehat{\sigma })$, showing in particular that they are
non-negative.

A {\it triple} consists of three boxes $(u,v,w)$ in $\augdg(\mu)$ such
that
\begin{equation}
\text{$w = d(u)$ and  $v\in \arm (u)$}.
\end{equation}
Given a filling $\widehat{\sigma }$ and boxes $x$, $y$ in $\augdg (\mu
)$, with $x<y$ in the reading order, set
\begin{equation}
\chi _{xy}(\widehat{\sigma }) = \begin{cases}
1&	\text{if $\widehat{\sigma }(x)>\widehat{\sigma }(y)$}\\
0&	\text{otherwise}.
\end{cases}
\end{equation}
If $(u,v,w)$ is a triple, then $v$ attacks $u$ and $w$, so we have
$\chi _{uv}(\widehat{\sigma })=1\Leftrightarrow (u,v)\in \Inv
(\widehat{\sigma })$, $\chi _{vw}(\widehat{\sigma })=1\Leftrightarrow
(v,w)\in \Inv (\widehat{\sigma })$, and $\chi_{uw} (\widehat{\sigma 
})=1\Leftrightarrow u\in \Des (\widehat{\sigma })$.  It follows from
the transitive law for inequalities that we always have
\begin{equation}
\chi _{uv}(\widehat{\sigma })+\chi _{vw} (\widehat{\sigma })- \chi
_{uw}(\widehat{\sigma })\in \{0,1 \}.
\end{equation}
We say that the triple $(u,v,w)$ is an {\it inversion triple} of
$\sigma$ if $\chi _{uv}(\widehat{\sigma })+\chi _{vw} (\widehat{\sigma
})- \chi _{uw}(\widehat{\sigma }) = 1$.  Otherwise $(u,v,w)$ is a {\it
co-inversion triple}.

\begin{lemma}\label{lem:attack-pairs}
Every pair of attacking boxes in $\augdg (\mu )$ occurs as either
$\{u,v \}$ or $\{v,w \}$ in a unique triple $(u,v,w)$, with the
exception that an attacking pair $\{(i,0),(i',0) \}$ in row $0$, such
that $i<i'$ and $\mu _{i}\leq \mu _{i'}$, does not belong to any
triple.
\end{lemma}

\begin{proof}
Consider an attacking pair $\{x=(i,j),y=(i',j) \}$ in the same row,
where $i<i'$.  If $\mu _{i}\leq \mu _{i'}$, then $x=(i,j)$ is in the
arm of $y=(i',j)$, so $(y,x,w=d(y))$ is a triple, provided $j\not =0$.
If $\mu _{i}>\mu _{i'}$, then $u=(i,j+1)$ is a box of $\skydg (\mu )$,
and $y$ is in the arm of $u$, so $(u,y,x)$ is a triple.  Moreover, in
the first case, $y$ is not in the arm of $u$ (possibly $u$ is not even
in the diagram), so $(u,y,x)$ is not a triple, while in the second
case, $x$ is not in the arm of $y$, so $(y,x,w=d(y))$ is not a
triple.  These are the only two sets of three boxes that might potentially
form a triple containing  $\{x,y \}$, so the triple with this property
is unique.

We leave the similar argument for an attacking pair in consecutive
rows to the reader.
\end{proof}

\begin{prop}\label{prop:inv}
Let $\sigma $ be a filling of $\mu $ (not necessarily non-attacking).
The number of inversion triples (resp.\ co-inversion triples) of
$\widehat{\sigma }$ is equal to $\inv (\widehat{\sigma })$ (resp.\
$\coinv (\widehat{\sigma })$).
\end{prop}

\begin{proof}
The number of inversion triples is given by the sum over all triples
\begin{equation}
\sum _{(u,v,w)}\bigl( \chi _{uv}(\widehat{\sigma })+\chi
_{vw}(\widehat{\sigma }) - \chi _{uw}(\widehat{\sigma })\bigr).
\end{equation}
Lemma~\ref{lem:attack-pairs} implies that the contribution to the sum
from the first two terms is equal to $|\Inv (\widehat{\sigma })|-
|\{i<j:\mu _{i}\leq \mu _{j} \}|$.  The contribution from the last
term is $-\sum _{u\in \Des (u)} a(u)$.  This proves that $\inv
(\widehat{\sigma })$ is equal to the number of inversion triples.
There is exactly one triple $(u,v,w)$ for every $u\in \skydg (\mu )$
and $v\in \arm (u)$.  Hence the total number of triples is equal to
$\sum _{u\in \skydg (\mu )} a(u)$, and it follows that $\coinv
(\widehat{\sigma })$ is equal to the number of co-inversion triples.
\end{proof}

A more pictorial characterization of inversion versus co-inversion
triples is sometimes useful.  Observe that the pattern of boxes in any
triple $(u,v,w)$ is one of the following two types:
\begin{equation}\label{e:triple-types}
\begin{array}{c@{\qquad \qquad }c}
\tableau{ u&       \\
          w& & &v    }\; , &  \tableau{ v& & &u\\
                                         & & &w   }\; , \\[2ex]
\text{Type I}  &              \text{Type II}
\end{array}
\end{equation}
with the further proviso that the column containing $u$, $w$ is
strictly taller than the column containing $v$ in Type I, and weakly
taller in Type II.  In either type, one easily checks the following
criterion.

\begin{lemma}\label{lem:orientation}
If $\sigma $ is a non-attacking filling, then a triple $(u,v,w)$ is a
co-inversion triple if and only if
$\widehat{\sigma}(u)<\widehat{\sigma}(v)<\widehat{\sigma}(w)$ or
$\widehat{\sigma}(v)<\widehat{\sigma}(w)<\widehat{\sigma}(u)$ or
$\widehat{\sigma}(w)<\widehat{\sigma}(u)<\widehat{\sigma}(v)$.
\end{lemma}

Informally speaking, a triple in a non-attacking filling is a
co-inversion triple if its entries increase {\it clockwise} in Type I,
or {\it counterclockwise} in Type II.  For example, the two
co-inversion triples in \eqref{e:example} are a Type I formed by the
$3$ and the $5$ in row $1$ with the $4$ in row $2$, and a Type II
formed by the $6$ and the $4$ in row $2$ with the $3$ in row $1$.

We may recast \eqref{e:main} as a formula for the ``opposite''
Macdonald polynomials $E_{\mu }(x;q^{-1},t^{-1})$, and interpret the
result conveniently in terms of inversion triples.  Define
\begin{equation}\label{e:anti-inv-maj}
\begin{aligned}
\maj '(\widehat{\sigma }) & = \bigl(\sum _{\substack{u\in \skydg (\mu)\\
\widehat{\sigma }(u)\not =\widehat{\sigma }(d(u))}} (l(u)+1) \bigr) -
\maj (\widehat{\sigma }) = \sum _{\substack{u\in \skydg (\mu)\\
\widehat{\sigma }(u) < \widehat{\sigma }(d(u))}} (l(u)+1), \\
\coinv '(\widehat{\sigma }) & =  \bigl(\sum _{\substack{u\in \skydg (\mu)\\
\widehat{\sigma }(u)\not =\widehat{\sigma }(d(u))}} a(u) \bigr) -
\coinv (\widehat{\sigma }).
\end{aligned}
\end{equation}
Then $\coinv '(\widehat{\sigma })$ is the number of inversion triples
with distinct entries---or, if one prefers, $\coinv '(\widehat{\sigma
}) = \coinv (\widehat{\sigma }')$ and $\maj '(\widehat{\sigma }) =
\maj (\widehat{\sigma }')$, where $\widehat{\sigma }'(u) =
n+1-\widehat{\sigma }(u)$.  In particular, $\coinv '(\widehat{\sigma
}) $ and $\maj '(\widehat{\sigma })$ are non-negative.  We have the
following corollary to Theorem~\ref{thm:main}.

\begin{cor}\label{cor:main-inverted}
\begin{equation}\label{e:main-inverted}
E_{\mu }(x;q^{-1},t^{-1}) = \sum _{\substack{\sigma \colon \mu
\rightarrow [n]\\
\text{non-attacking}}} x^{\sigma } q^{\maj '(\widehat{\sigma })}
t^{\coinv '
(\widehat{\sigma })} \prod _{\substack{u\in \skydg (\mu )\\
\widehat{\sigma }(u)\not =\widehat{\sigma }(d(u))}}
\frac{1-t}{1-q^{l(u)+1}t^{a(u)+1}}.
\end{equation}
\end{cor}

\section{Proof of the formula}
\label{sec:proof}

This section is devoted to the proof of Theorem~\ref{thm:main}.  To
avoid confusion between the two quantities which we want to prove are
equal, we will henceforth denote the right-hand side of \eqref{e:main}
by $C_{\mu }(x;q,t)$.  We will prove the theorem by verifying that the
special cases of \eqref{e:intertwiners} and \eqref{e:Psi} isolated in
Lemma~\ref{lem:Knop-recurrence} hold with $C_{\mu }$ in place of
$E_{\mu }$.

\subsection{The cyclic shift symmetry $\pi $.}
\label{subsec:cyclic-shift}

Our first task is to verify that $C_{\mu }$ satisfies \eqref{e:Psi}.
It will be convenient to extend the symmetry $\pi $ in
\eqref{e:pi-Psi} first to {\it diagrams}, defining
\begin{equation}
\pi \colon \augdg (\mu )\rightarrow \augdg (\pi (\mu ))
\end{equation}
by $\pi ((i,j)) = (i+1,j)$ for $i<n$, and $\pi ((n,j)) = (1,j+1)$,
also to {\it values} of fillings, defining 
\begin{equation}
\pi \colon [n]\rightarrow [n]
\end{equation}
by $\pi (i) = i+1$ for $i<n$, and $\pi (n) = 1$, and finally to {\it
fillings}, defining $^{\pi }\sigma $ to be the unique filling of
$\skydg (\pi (\mu ))$ whose augmented filling satisfies
\begin{equation}
\widehat{^{\pi }\sigma }(\pi (u)) = \pi (\widehat{\sigma }(u))
\end{equation}
for all $u\in \augdg (\mu )$.  Note that this definition is compatible
with the requirement that $\widehat{^{\pi }\sigma}((j,0)) = j$, so
$\widehat{^{\pi }\sigma }$ is a well-defined augmented filling.
Also note that $\pi \colon \augdg (\mu )\rightarrow \augdg (\pi (\mu
))$ maps attacking pairs to attacking pairs; hence if $\sigma $ is
non-attacking, then so is $^{\pi }\sigma $.

\begin{lemma}\label{lem:shift-arm}
Given $u\in \skydg (\mu )$, $v\in \augdg (\mu )$, we have $v\in \arm
(u)$ if and only if $\pi (v)\in \arm (\pi (u))$.
\end{lemma}

\begin{proof}
Immediate from the definitions.
\end{proof}

\begin{prop}\label{prop:shift-symmetry}
Let $\sigma \colon \skydg (\mu )\rightarrow [n]$ be a non-attacking
filling.  Then $\coinv (\widehat{^{\pi }\sigma })=\coinv
(\widehat{\sigma })$, and $\maj (\widehat{^{\pi }\sigma }) = \maj
(\widehat{\sigma })+\mu _{n}-r$, where $r = |\sigma ^{-1}(\{n \})|$.
\end{prop}

\begin{proof}
By Lemma~\ref{lem:shift-arm}, if $(u,v,w)$ is a triple in $\augdg (\mu
)$, then $(u',v',w') = (\pi (u),\pi (v),\pi (w))$ is a triple in
$\augdg (\pi (\mu ))$.  From Lemma~\ref{lem:orientation}, we see that
$(u',v',w')$ is a co-inversion triple of $\widehat{^{\pi }\sigma }$ if
and only if $(u,v,w)$ is a co-inversion triple of $\widehat{\sigma }$.
Moreover, every triple in $ \augdg (\pi (\mu ) )$ is the image under
$\pi $ of a triple in $\augdg (\mu )$, except for those triples
$(u',v',w')$ in which $w'=(1,0)$.  But these exceptional triples
satisfy $\widehat{^{\pi }\sigma }(u') = \widehat{^{\pi }\sigma
}(w')=1$, so they are not co-inversion triples.  This proves that
$\coinv (\widehat{^{\pi }\sigma })=\coinv (\widehat{\sigma })$.

For the major index, let $S = \sigma ^{-1}(\{n \}) = \{u\in \skydg
(\mu ):\widehat{\sigma }(u) = n \}$, and let $S' = \{u\in \skydg (\mu
):\widehat{\sigma } (d(u)) = n \}$.  For $u\in \skydg (\mu )$, we have
$\pi (u)\in \Des (\widehat{^{\pi }\sigma })$ if and only if either
$u\in \Des (\widehat{\sigma })\setminus S$, or $u\in S'\setminus S$.
Note that the only box of $\skydg (\pi (\mu ))$ not in the image
$\pi(\skydg (\mu ))$ is $(1,1)$, which is not a descent of
$\widehat{^{\pi }\sigma }$, and that $l(\pi (u)) = l(u)$ for all $u\in
\skydg (\mu )$.  Therefore,
\begin{equation}\label{e:maj-pi-sigma}
\maj (\widehat{^{\pi }\sigma }) = \sum _{u\in \Des (\widehat{\sigma
})\setminus S} (l(u)+1) +\sum _{u\in S'\setminus S} (l(u)+1).
\end{equation}
Now, $S\cap \Des (\widehat{\sigma }) = S\setminus S'$, and using this,
\eqref{e:maj-pi-sigma} is equivalent to
\begin{equation}\label{e:maj-pi-sigma-II}
\maj (\widehat{^{\pi }\sigma }) = \sum _{u\in \Des (\widehat{\sigma
})}(l(u)+1) +\sum _{u\in S'}(l(u)+1) - \sum _{u\in S} (l(u)+1).
\end{equation}
Observe that if $\mu _{n}>0$, then $S'$ contains the box $v=(n,1)$,
and that $l(v)+1=\mu _{n}$.  The map $u'\rightarrow d(u')$ is a
bijection from $S'\setminus \{v \}$ onto the set of boxes $u\in S$
such that $l(u)\not =0$.  Hence the second term in
\eqref{e:maj-pi-sigma-II} is equal to
\begin{equation}
\mu _{n}+\sum _{\substack{u\in
S\\ l(u)\not =0}} l(u) = \mu _{n}+\sum _{u\in S}l(u),
\end{equation}
so we obtain $\maj (\widehat{^{\pi }\sigma }) = \maj (\widehat{\sigma
})+\mu _{n}-|S|$, as desired.
\end{proof}

\begin{cor}\label{cor:shift-symmetry}
The combinatorial expression $C_{\mu }$ satisfies \eqref{e:Psi},
that is,
\begin{equation}\label{e:Psi-for-C}
C_{\pi (\mu )}(x;q,t) = q^{\mu _{n}}\Psi C_{\mu }(x;q,t) \defeq q^{\mu
_{n}} x_{1} C_{\mu }(x_{2},\ldots,x_{n},q^{-1}x_{1};q,t).
\end{equation}
\end{cor}

\begin{proof}
Every non-attacking filling of $\pi (\mu )$ satisfies $\sigma
((1,1))=1$, since $(1,1)$ attacks all the boxes $(j,0)$ for
$j=2,\ldots,n$.  Therefore $\pi $ is a bijection from non-attacking
fillings of $\mu $ to non-attacking fillings of $\pi (\mu )$.
Clearly, $\Psi (x^{\sigma }) = q^{-r}x^{(^{\pi }\sigma )}$, where $r =
|\sigma ^{-1}(\{n \})|$.  Hence $q^{\mu _{n}}\Psi (x^{\sigma })q^{\maj
(\widehat{\sigma })}t^{\coinv (\widehat{\sigma })} = x^{(^{\pi }\sigma
)}q^{\maj (\widehat{^{\pi }\sigma })}t^{\coinv (\widehat{^{\pi }\sigma
})}$, by Proposition~\ref{prop:shift-symmetry}.  Moreover, $\pi $
induces a bijection from $\{u\in \skydg (\mu ):\widehat{\sigma
}(u)\not =\widehat{\sigma }(d(u)) \}$ to $\{u\in \skydg (\pi (\mu
)):\widehat{^{\pi }\sigma } (u)\not =\widehat{^{\pi }\sigma } (d(u))
\}$.  It is obvious that $l(\pi (u))=l(u)$, and
Lemma~\ref{lem:shift-arm} gives $a(\pi (u))=a(u)$.  Hence the term of
$C_{\pi (\mu )}$ corresponding to $^{\pi }\sigma $ contains the same
factor
\begin{equation}
\prod _{\substack{u\in \skydg (\mu )\\
\widehat{\sigma }(u)\not =\widehat{\sigma }(d(u))}}
\frac{1-t}{1-q^{l(u)+1}t^{a(u)+1}}
\end{equation}
as the term of $C_{\mu }$ corresponding to $\sigma $.  This shows that
\eqref{e:Psi-for-C} holds term by term.
\end{proof}

\subsection{Partial symmetry of generating functions for fillings}
\label{subsec:symmetry}

It remains to verify that $C_{\mu }$ satisfies suitable cases of
\eqref{e:intertwiners}.  We cannot do this simply by applying the
operator $T_{i}$ in \eqref{e:cherednik} to the formula for $C_{\mu }$,
as the resulting expressions are intractable.  Instead, we will take
an indirect approach, based on the combinatorial machinery associated
with LLT polynomials, which we used in \cite{HaHaLo05} to prove that
the formulas given there for symmetric Macdonald polynomials are
indeed symmetric.  Clearly we must expect less in the non-symmetric
case, but it turns out that we will be able to decompose $C_{\mu }$
into parts that are symmetric, or nearly so, in the two variables
$x_{i}$ and $x_{i+1}$, which is enough to compute $T_{i}C_{\mu }$.

Fix the signed alphabet
\begin{equation}
\Acal (n) = \{1 < \overline{1} <2 <\overline{2} < \cdots < n < \overline{n}\}
\end{equation}
with the ordering that we denoted by $<_{1}$ in \cite{HaHaLo05}.  Let
$S\subseteq \NN ^{2}$ be an arbitrary finite lattice-square diagram.
A {\it signed filling} is a map $\sigma:S\rightarrow \Acal (n)$.
Descents and inversions of $\sigma $ are defined as in \S
\ref{subsec:des-maj}--\ref{subsec:inv}, with the modification that
pairs of equal negative entries count as descents or inversions.
Denote
\begin{equation}
\begin{aligned}
\Inv (\sigma ) & = \{\text{inversions of $\sigma $} \} \\
\Des (\sigma ) & = \{\text{descents of $\sigma $} \}.
\end{aligned}
\end{equation}

In \cite{HaHaLo05} we proved the following basic results. (The
definitions of reading order and attacking pairs in \cite{HaHaLo05}
differ from those used here by a change of coordinates, and the
results are stated for partition diagrams.  The coordinate change is
inconsequential, however, and the proofs also go through for arbitrary
diagrams.)

\begin{lemma}[{\cite[Proposition 3.4]{HaHaLo05}}]\label{lem:LLT}
Given a lattice-square diagram $S$ and a subset $D\subseteq S$, the
polynomial 
\begin{equation}
F_{S,D}(x;t) \defeq \sum _{\substack{\sigma \colon S\rightarrow [n]\\
\Des (\sigma ) = D}} x^{\sigma } t^{| \Inv (\sigma ) |}
\end{equation}
is a symmetric function; in fact it is an LLT polynomial.
\end{lemma}

\begin{lemma}[see proof of {\cite[eqs. (38), (40)]{HaHaLo05}}]
\label{lem:LLT-transformed} For a signed filling $\sigma \colon
S\rightarrow \Acal (n)$, let $p(\sigma )$ and $m(\sigma )$ denote the
number of positive and negative entries $\sigma (u)$, respectively.
Then with $S$ and $D$ as in Lemma~\ref{lem:LLT}, and with $X =
x_{1}+\cdots +x_{n}$, we have, using plethystic notation,
\begin{equation}
F_{S,D}[X(t-1);t] = \sum _{\substack{\sigma \colon S\rightarrow \Acal (n)\\
\Des (\sigma ) = D}} x^{|\sigma |}(-1)^{m(\sigma )} t^{p(\sigma 
)+| \Inv(\sigma ) |}.
\end{equation}
In particular, the sum on the right-hand side is a symmetric function.
\end{lemma}

Define a signed filling $\sigma $ to be {\it non-attacking} if
$|\sigma |$ is non-attacking.

\begin{lemma}[{\cite[Lemma 5.1 and remark after it]{HaHaLo05}}]
\label{lem:LLT-non-attacking}
The polynomial $F_{S,D}[X(t-1);t]$ in Lemma~\ref{lem:LLT-transformed} is
also given by
\begin{equation}\label{e:non-attacking-form}
F_{S,D}[X(t-1);t] = \sum _{\substack{\sigma \colon S\rightarrow \Acal (n)\\
\Des (\sigma ) = D\\
\text{$\sigma $ non-attacking}}} x^{|\sigma |}(-1)^{m(\sigma )}
t^{p(\sigma )+| \Inv (\sigma ) |}.
\end{equation}
In particular, the sum on the right-hand side is a symmetric function.
\end{lemma}

Now assign arbitrary ``arm'' and ``leg'' values $a(u)$, $l(u)$ to each
$u\in S$ such that $d(u)\in S$, and
define for any filling $\sigma $ of $S$,
\begin{equation}\label{e:generalized-maj-inv}
\begin{gathered}
\inv (\sigma ) = |\Inv (\sigma )| -\sum _{u\in \Des (\sigma )} a(u)\\
\maj (\sigma ) = \sum _{u\in \Des (\sigma )} (l(u)+1).
\end{gathered}
\end{equation}

\begin{lemma}\label{lem:maj-inv-general}
We have the identity
\begin{multline}\label{e:maj-inv-general}
\sum _{\substack{\sigma \colon S\rightarrow \Acal (n)\\
\text{non-attacking}}} x^{|\sigma |}(-1)^{m(\sigma )}q^{\maj (\sigma
)}t^{|S|-(p(\sigma)+\inv (\sigma ))}\\
 =  \sum _{\substack{\sigma \colon S\rightarrow [n]\\
\text{non-attacking}}} x^{\sigma }q^{\maj (\sigma )}t^{-\inv (\sigma
)} \prod _{\substack{u,d(u)\in S\\
\sigma (u)=\sigma (d(u))}} \left(1-q^{l(u)+1}t^{a(u)+1} \right) \prod
_{\substack{u\in S\\ \sigma (u)\not =\sigma (d(u))}} (1-t),
\end{multline}
where boxes $u\in S$ such that $d(u)\not \in S$ contribute to the last
factor.  In particular, the sum on the right-hand side is a symmetric
function.
\end{lemma}

\begin{proof}
On the left-hand side, consider the partial sum over terms for which
$|\sigma |$ is a given filling $\tau \colon S\rightarrow [n]$.  Set $V
= \{u\in S: d(u)\in S,\, \tau (u)=\tau (d(u)) \}$.  Since $\tau $ is
non-attacking, we have $\Inv (\sigma ) = \Inv (\tau )$.  We also have
$\Des (\sigma ) = \Des (\tau )\cup \{u\in V: \sigma (u)\in
\overline{[n]} \}$.  In particular, when $\sigma =\tau $, we have
$\Des (\sigma ) = \Des (\tau )$, and the corresponding term on the
left-hand side of \eqref{e:maj-inv-general} is $x^{\tau }q^{\maj (\tau
)}t^{-\inv (\tau )}$.  In general, this term gets multiplied by a
factor $q^{l(u)+1}t^{a(u)}\cdot (-t)$ for each $u\in V$ such that
$\sigma (u)$ is negative, and by a factor $-t$ for each $u\in
S\setminus V$ such that $\sigma (u)$ is negative.  Summing over all
choices of signs gives the term
\begin{equation}
x^{\tau }q^{\maj (\tau )}t^{-\inv (\tau )} \prod _{u\in V}
\left(1-q^{l(u)+1}t^{a(u)+1} \right) \prod _{u\in S\setminus V} (1-t),
\end{equation}
and summing this over all $\tau $ gives \eqref{e:maj-inv-general}.  On
the left-hand side of \eqref{e:maj-inv-general}, the partial sum for
each fixed descent set $\Des (\sigma ) = D$ is a constant multiple of
\eqref{e:non-attacking-form} (with $t^{-1}$ in place of $t$).  Hence
each side of \eqref{e:maj-inv-general} is a symmetric function.
\end{proof}

We need to work in the more general setting of augmented fillings
$\widehat{\sigma }$ that take prescribed values on the extra boxes in
a larger diagram $\widehat{S}\supseteq S$.  In this setting, of
course, we lose the full symmetry, but we sometimes keep a partial
symmetry. 

\begin{prop}\label{prop:augmented-symmetry}
Suppose given two disjoint lattice-square diagrams $S$ and $T$, two
disjoint subsets $Y,Z\in S$, and a filling $\tau \colon T\rightarrow
[n]$ such that $\tau (T)\cap \{i, i+1\} = \emptyset$.  Let
$\widehat{S} = S\cup T$ and fix arm and leg values $a(u)$, $l(u)$ for
each box $u\in S$ such that $d(u)\in \widehat{S}$.  For any filling
$\sigma \colon S\rightarrow [n]$, set $\widehat{\sigma } = \sigma \cup
\tau $.  With $\maj (\widehat{\sigma })$ and $\inv (\widehat{\sigma 
})$ defined as in \eqref{e:generalized-maj-inv}, the sum
\begin{equation}\label{e:augmented-symmetry}
\sum _{\substack{\sigma \colon S\rightarrow [n]\\
\text{$\widehat{\sigma }$ non-attacking}\\
\sigma (Y)\cap\{i, i+1\} = \emptyset \\
\sigma (Z)\subseteq \{i, i+1\} }} x^{\sigma }q^{\maj (\widehat{\sigma 
})} t^{-\inv (\widehat{\sigma })}
\prod _{\substack{u\in S,d(u)\in \widehat{S}\\
\widehat{\sigma }(u)=\widehat{\sigma }(d(u))}}
\left(1-q^{l(u)+1}t^{a(u)+1} \right) \prod _{\substack{u\in
S\\ \widehat{\sigma }(u)\not =\widehat{\sigma }(d(u))}}
(1-t),
\end{equation}
with any boxes $u\in S$ such that $d(u)\not \in \widehat{S}$ included
in the last factor, is symmetric in $x_{i}$ and $x_{i+1}$.
\end{prop}

\begin{proof}
Consider the partial sum over fillings for which the set $Z' = \sigma
^{-1}(\{i, i+1 \})$ and the restriction of $\sigma $ to $S\setminus
Z'$ are fixed.  To evaluate this sum, we may as well adjoin $Z'$ to
$Z$ and $S\setminus Z'$ to $T$, thus reducing the problem to the case
that $Y = \emptyset $, $Z = S$, and the sum is over fillings $\sigma
\colon S \rightarrow \{i, i+1\}$.

For each $u\in T$, the entry $\tau (u)$ stands in the same relative
order to both $i$ and $i+1$.  Hence the only descents, inversions, and
pairs $\widehat{\sigma }(u) = \widehat{\sigma }(d(u))$ which depend on
$\sigma $ are the ones created by pairs of boxes that are both in
$S$.  The whole sum \eqref{e:augmented-symmetry} is therefore the
product of a fixed polynomial $c(x;q,t)$ and the corresponding sum for
the diagram $S$ without the augmentation.  In other words, the problem
is further reduced to the case $T=\emptyset $, which is
Lemma~\ref{lem:maj-inv-general}.
\end{proof}

\begin{cor}\label{cor:augmented-symmetry}
Keep the notation of Proposition~\ref{prop:augmented-symmetry}.  The sum
\begin{equation}\label{e:augmented-symmetry-II}
\sum _{\substack{\sigma \colon S\rightarrow [n]\\
\text{$\widehat{\sigma }$ non-attacking}\\
\sigma (Y)\cap\{i, i+1\} = \emptyset \\
\sigma (Z)\subseteq \{i, i+1\} }} x^{\sigma }q^{\maj (\widehat{\sigma 
})} t^{-\inv (\widehat{\sigma })} \prod _{\substack{u\in S, d(u)\in
\widehat{S}\\ \widehat{\sigma }(u)\not =\widehat{\sigma }(d(u))}}
\frac{(1-t)}{(1-q^{l(u)+1}t^{a(u)+1})},
\end{equation}
is symmetric in $x_{i},x_{i+1}$.
\end{cor}

\begin{proof}
Divide \eqref{e:augmented-symmetry} by $\prod _{u\in S,d(u)\in
\widehat{S}}\left(1-q^{l(u)+1}t^{a(u)+1} \right)\prod _{u\in
S,d(u)\not \in S}(1-t)$.
\end{proof}

\subsection{Operators $T_{i}$ and symmetry}
\label{subsec:hecke}

The following lemma will enable us to apply the machinery in \S
\ref{subsec:symmetry}.

\begin{lemma}\label{lem:magic}
For any $G_{1},G_{2}\in \QQ (q,t)X$, and $0<i<n$, the following
conditions are equivalent:
\begin{itemize}
\item [(i)] $G_{2} = T_{i}G_{1}$;
\item [(ii)] $G_{1}+G_{2}$ and $tx_{i+1}G_{1}+x_{i}G_{2}$ are
symmetric in $x_{i}, x_{i+1}$.
\end{itemize}
\end{lemma}

\begin{proof}
Let $A = (\QQ (q,t)X)^{s_{i}}$ be the subring of Laurent polynomials
which are symmetric in $x_{i},x_{i+1}$.  For every $f\in \QQ(q,t)X$,
there exist $a,b\in A$ such that
\begin{equation} 
f = a +  x_{i} b,\quad a,b\in A.
\end{equation}
Specifically, we can take $b=(f-s_{i}(f))/(x_{i}-x_{i+1})$ and
$a=(x_{i+1} f-x_{i}s_{i}(f))/(x_{i+1}-x_{i})$.  Now, $T_{i}$ is an
$A$-linear operator, satisfying the identity
\begin{equation}
T_{i}(a+x_{i}b) = ta+x_{i+1}b,\quad a,b\in A,
\end{equation}
since $T_{i}(1) = t$ and $T_{i}(x_{i})= x_{i+1}$.

Suppose (i) holds.  Express $G_{1}$ in the form
\begin{equation}\label{e:G1}
G_{1} = a + x_{i}b.
\end{equation}
Then we have
\begin{align}
G_{2} & = t a + x_{i+1} b, \label{e:G2}\\
G_{1}+G_{2} & = (1+t) a + (x_{i}+x_{i+1}) b, \label{e:G-sym-I}\\
t x_{i+1}G_{1} + x_{i} G_{2} & = t(x_{i}+x_{i+1}) a +
(1+t)x_{i}x_{i+1} b. \label{e:G-sym-II}
\end{align}
The right-hand sides of \eqref{e:G-sym-I}--\eqref{e:G-sym-II} are
manifestly symmetric, giving (ii).  Conversely, suppose (ii) holds.
We can view \eqref{e:G-sym-I}--\eqref{e:G-sym-II} as a non-singular
system of linear equations for unknown rational functions $a$ and $b$,
with symmetric coefficients.  Hence there exist $a,b\in \QQ
(q,t,x)^{s_{i}}$ such that \eqref{e:G-sym-I}--\eqref{e:G-sym-II} hold,
and therefore \eqref{e:G1}--\eqref{e:G2} hold also.  The operator
$T_{i}$ extends to $\QQ (q,t,x)$ and is linear over $\QQ
(q,t,x)^{s_{i}}$, so it follows that $T_{i}G_{1} = ta+ x_{i+1} b =
G_{2}$.
\end{proof}

\subsection{Conclusion of the proof}
\label{subsec:proof}

Fix $\mu $ and $i$ such that $\mu _{i+1}=0$, as in
Lemma~\ref{lem:Knop-recurrence}(b).  We must show that $C_{\mu }$
satisfies \eqref{e:intertwiners}, or equivalently,
\eqref{e:arm-leg-intertwiners}.  Explicitly, let $u = (i,1)\in \skydg
(\mu )$.  We are to prove
\begin{equation}\label{e:C-mu-intertwiner}
C_{s_{i}(\mu )}(x;q,t) = \left(T_{i} + \frac{1-t}{1-q^{\mu
_{i}}t^{a(u)}} \right) C_{\mu }(x;q,t).
\end{equation}
Let 
\begin{equation}
C_{\mu } = G_{0}+G_{1},
\end{equation}
where $G_{0}$ is the partial sum in \eqref{e:main} over fillings
$\sigma $ such that $\sigma (u)\not =i$, and $G_{1}$ is the sum over
fillings such that $\sigma (u) = i$.

\begin{lemma}\label{lem:symm-I}
$G_{0}$ is symmetric in $x_{i}$, $x_{i+1}$.
\end{lemma}

\begin{proof}
In Corollary~\ref{cor:augmented-symmetry}, take $S=\skydg (\mu )$, $T
= \{(j,0):j\in [n]\setminus \{i,i+1 \} \}$, $\tau ((j,0)) = j$, $Y =
\{(j,1)\in \skydg (\mu ):j\leq i \}$, and $Z = \emptyset $.  The only
conditions imposed on a non-attacking filling $\sigma $ of $\mu $ by
the augmentation in boxes $(i,0)$ and $(i+1,0)$ are that $\sigma 
(y)\not \in \{i,i+1 \}$ for $y\in Y\setminus \{u \}$ and $\sigma 
(u)\not =i+1$.  Hence the non-attacking augmented fillings
$\widehat{\sigma }'$ of $\widehat{S} = \augdg (\mu )\setminus
\{(i,0),(i+1,0) \}$ which satisfy $\widehat{\sigma }'(y)\not \in
\{i,i+1 \}$ for $y\in Y$ are precisely the restrictions to
$\widehat{S}$ of the augmentations $\widehat{\sigma }$ of fillings
$\sigma $ of $\mu $ such that $\sigma (u)\not =i$.

Now, $G_{0}$ is the partial sum in \eqref{e:main} over the above
fillings $\sigma $, while the expression in
\eqref{e:augmented-symmetry-II} is a sum over their restrictions
$\widehat{\sigma }'$.  Clearly, $\maj (\widehat{\sigma }') = \maj
(\widehat{\sigma })$, since $u$ is never a descent of $\widehat{\sigma
}$.  For every box $(j,1)$ in row $1$ of $\skydg (\mu )$, we have
$\sigma ((j,1))\leq j$, since $(j,1)$ attacks $(k,0)$ for all $k>j$.
Hence $\widehat{\sigma }$ has no inversions between boxes in row $1$
and boxes in row $0$, and therefore $\Inv (\widehat{\sigma })$ is the
union of $\Inv (\widehat{\sigma }')$ with the fixed set of $2n-3$
additional inversions in row $0$ involving the boxes $(i,0)$ or
$(i+1,0)$.  This shows that $\coinv (\widehat{\sigma })$ differs from
$-\inv (\widehat{\sigma }')$ by a constant independent of $\sigma $.
Finally, since $\sigma (u)\not =i$, each term of $G_{0}$ contains a
factor $(1-t)/(1-q^{\mu _{i}}t^{a(u)})$ not in the corresponding term
of \eqref{e:augmented-symmetry-II}, but the remaining factors are the
same for $\widehat{\sigma }$ and $\widehat{\sigma }'$.  Hence $G_{0}$
is a constant multiple of the expression in
\eqref{e:augmented-symmetry-II}.
\end{proof}

Now define $s_{i}\colon \skydg (\mu )\rightarrow \skydg (s_{i}(\mu ))$
by $s_{i}((j,k)) = (s_{i}(j),k)$.  Let $G_{2}$ be the sum of the terms
of $C_{s_{i}(\mu )}$ corresponding to fillings $\sigma $ such that
$\sigma (s_{i}(u)) = i+1$.

\begin{lemma}\label{lem:symm-II}
$G_{1}+G_{2}$ is symmetric in $x_{i}$, $x_{i+1}$.
\end{lemma}

\begin{proof}
In Corollary~\ref{cor:augmented-symmetry}, take $S=\skydg (\mu )$, $T
= \{(j,0):j\in [n]\setminus \{i,i+1 \} \}$, $\tau ((j,0)) = j$, $Y =
\{(j,1)\in \skydg (\mu ):j< i \}$, and $Z =\{u \}$.  By similar
reasoning to that in the proof of the previous lemma, the fillings
$\sigma '\colon S\rightarrow [n]$ in \eqref{e:augmented-symmetry-II}
such that $\sigma '(u)=i+1$ correspond bijectively to non-attacking
fillings $\sigma = \sigma '\circ s_{i}$ of $s_{i}(\mu )$ such that
$\sigma (s_{i}(u))=i+1$.  Again, $\coinv (\widehat{\sigma })$ differs
from $-\inv (\widehat{\sigma }')$ by a constant $C$ independent of
$\sigma $.  Moreover, $\widehat{\sigma }(d(s_{i}(u))) = i+1 =
\widehat{\sigma }(s_{i}(u))$, and for every $v\in \skydg (\mu
)\setminus \{u \}$, we have $l(s_{i}(v))=l(v)$, $a(s_{i}(v))=a(v)$.
Hence the terms in \eqref{e:augmented-symmetry-II} for $\sigma
'(u)=i+1$ sum to $t^{-C}G_{2}$.  Similarly, the terms for $\sigma
'(u)=i$ sum to $t^{-C}G_{1}$, so the expression in
\eqref{e:augmented-symmetry-II} is equal to $t^{-C}(G_{1}+G_{2})$.
\end{proof}

\begin{lemma}\label{lem:symm-III}
$tx_{i+1}G_{1}+x_{i}G_{2}$ is symmetric in $x_{i}$, $x_{i+1}$.
\end{lemma}

\begin{proof}
Let $v = (i+1,0)$.  In Corollary~\ref{cor:augmented-symmetry}, take
$S=\skydg (\mu )\cup \{v \}$, $T = \{(j,0):j\in [n]\setminus \{i,i+1
\} \}$, $\tau ((j,0)) = j$, $Y = \{(j,1)\in \skydg (\mu ):j< i \}$,
and $Z =\{u, v \}$.  Every non-attacking filling $\sigma '\colon
S\rightarrow [n]$ such that $\sigma '(v)=i$ has $\sigma '(u)=i+1$.  As
in the proof of the previous lemma, these fillings correspond
bijectively to fillings $\sigma $ of $s_{i}(\mu )$ such that $\sigma
(s_{i}(u)) = i+1$, and there is a constant $C$ such that $\coinv
(\widehat{\sigma }) = C-\inv (\widehat{\sigma }')$.  Since $\sigma
'(v)=i$, we have $x^{\sigma '} = x_{i}x^{\sigma }$, so the
contribution to \eqref{e:augmented-symmetry-II} from these terms is
$t^{-C}x_{i}G_{2}$.  Note that $C$ includes a contribution from the
inversion $(u,v)\in \Inv (\widehat{\sigma }')\setminus \Inv
(\widehat{\sigma })$.

The remaining terms, with $\sigma '(v)=i+1$ and $\sigma '(u)=i$,
correspond to fillings $\sigma $ of $\mu $ such that $\sigma (u)=i$.
In this case, $(u,v)\not \in \Inv (\widehat{\sigma })$, so we now have
$\coinv (\widehat{\sigma }) = C-1-\inv (\widehat{\sigma }')$.  The
contribution from these terms is therefore $t^{1-C}x_{i+1}G_{1}$, and
the expression in \eqref{e:augmented-symmetry-II} is equal to
$t^{-C}(tx_{i+1}G_{1}+x_{i}G_{2})$.
\end{proof}

To complete the proof of Theorem~\ref{thm:main}, observe that
$T_{i}G_{0} = tG_{0}$, by Lemma~\ref{lem:symm-I}, and $T_{i}G_{1} =
G_{2}$, by Lemmas~\ref{lem:magic}, \ref{lem:symm-II} and
\ref{lem:symm-III}.  This yields
\begin{equation}\label{e:operator-on-C-mu}
\left(T_{i} + \frac{1-t}{1-q^{\mu _{i}}t^{a(u)}} \right) C_{\mu } =
\frac{1-q^{\mu _{i}}t^{a(u)+1}}{1-q^{\mu _{i}}t^{a(u)}} G_{0} +
\frac{1-t}{1-q^{\mu _{i}}t^{a(u)}} G_{1} + G_{2}.
\end{equation}
The fillings $\sigma '$ of $s_{i}(\mu )$ such that $\sigma 
'(s_{i}(u))\not \in \{i,i+1 \}$ correspond bijectively to fillings
$\sigma $ of $\mu $ such that $\sigma (u)\not =i$, by the rule $\sigma 
'=\sigma \circ s_{i}$.  Since $a(s_{i}(u)) = a(u)-1$, and $a(s_{i}(v))
= a(v)$ for all $v\not =u$ in $\skydg (\mu )$, the term in
$C_{s_{i}(\mu )}$ corresponding to $\sigma '$ differs from the term in
$C_{\mu }$ corresponding to $\sigma $ by the factor
\begin{equation}
\frac{1-q^{\mu _{i}}t^{a(u)+1}}{1-q^{\mu _{i}}t^{a(s_{i}(u))+1}} =
\frac{1-q^{\mu _{i}}t^{a(u)+1}}{1-q^{\mu _{i}}t^{a(u)}}.
\end{equation}
The first term in \eqref{e:operator-on-C-mu} is the sum of these
terms.  Similarly, the second term in \eqref{e:operator-on-C-mu} is
the sum of terms of $C_{s_{i}(\mu )}$ corresponding to $\sigma '$ such
that $\sigma '(s_{i}(u)) = i$.  The third term in
\eqref{e:operator-on-C-mu} is by definition the sum of the remaining
terms in $C_{s_{i}(\mu )}$, corresponding to $\sigma '$ such that
$\sigma '(s_{i}(u)) = i+1$.  Hence \eqref{e:operator-on-C-mu} gives
\eqref{e:C-mu-intertwiner}.

\section{Comparison with symmetric Macdonald polynomials}
\label{sec:comparison}

In this section we indicate some of the connections between the
combinatorial theory of non-symmetric Macdonald polynomials and the
corresponding theory in the symmetric case.  We will refer to three
variants of symmetric Macdonald polynomials: the monic forms
$P_{\lambda }(x;q,t)$ and integral forms $J_{\lambda }(x;q,t)$, as
defined in \cite{Mac95}, and the transformed integral forms
$\tilde{H}_{\lambda }(x;q,t)$ as in \cite{HaHaLo05}, \cite{HaHaLo05a}.
Throughout, $\lambda $ denotes a partition $(\lambda _{1}\geq \cdots
\geq \lambda _{n})$, possibly with some parts equal to zero.

\subsection{More general formulas for $\tilde{H}_{\lambda }$}
\label{subsec:more-H-formulas}

The main theorem in our previous paper \cite[Theorem 2.2]{HaHaLo05}
was a combinatorial formula for $\tilde{H}_{\lambda }$, expressed as a
sum over fillings of the diagram of $\lambda $.  By virtue of the
identity
\begin{equation}
\tilde{H}_{\lambda '}(x;q,t) = \tilde{H}_{\lambda }(x;t,q),
\end{equation}
the same thing can also be expressed as a sum over fillings of the
diagram of the transpose partition $\lambda '$, that is, of $\skydg
(\lambda )$.  We shall now generalize this to a sum over fillings of
any rearrangement of the columns of $\skydg (\lambda )$.  Our original
result is equivalent to the special case of the following theorem in
which $\mu $ is weakly increasing, {\it i.e.}, $\mu $ is the
anti-dominant weight in the Weyl group orbit of $\lambda $.

\begin{thm}\label{thm:general-H-formula}
For any composition $\mu \in \NN ^{n}$ and any $m>0$, define
\begin{equation}\label{e:D-mu}
D_{\mu }(x_{1},\ldots,x_{m};q,t) = \sum _{\sigma \colon \skydg (\mu
)\rightarrow [m]} x^{\sigma }q^{\maj (\sigma )}t^{\inv (\sigma )},
\end{equation}
where the sum is over {\it unaugmented, possibly attacking} fillings
$\sigma $, and $\maj (\sigma )$, $\inv (\sigma )$ are as in
\eqref{e:generalized-maj-inv}, using the arm and leg values defined by
\eqref{e:arms}--\eqref{e:a(u)-l(u)}.  (Note that, although the arm of
a box in row $1$ is defined with reference to the augmented diagram,
these arms play no role in formula \eqref{e:D-mu}, since a box in row
$1$ is never a descent of the unaugmented filling $\sigma $.)  Then
\begin{equation}\label{e:D-mu=H-lambda}
\tilde{H}_{\lambda }(x;q,t) = D_{\mu }(x;q,t),
\end{equation}
where $\mu $ is any rearrangement of the parts of the partition
$\lambda $.
\end{thm}

\begin{proof}
As a matter of fact, the proof of \cite[Theorem 2.2]{HaHaLo05} goes
through essentially unchanged in this more general setting.  We will
just briefly point out what properties of fillings of a partition
diagram were actually used in \cite{HaHaLo05}, in order to see that
the arguments there remain valid.  The reader should bear in mind that
for compatibility with the non-symmetric case, in \eqref{e:D-mu} we
use the right-to-left reading word, instead of the left-to-right
reading word used in \cite{HaHaLo05}.

First, the proof of \cite[Theorem~3.1]{HaHaLo05} applies verbatim to
show that $D_{\mu }$ is a symmetric function, and the derivation of
\cite[eqs.~(40)--(41)]{HaHaLo05} applies equally to $D_{\mu }$.

Next, in \cite[\S 5.1]{HaHaLo05}, we used two properties of the
diagram and the statistics $\inv (\sigma )$ and $\maj (\sigma )$: (i)
the fact that $\inv (\sigma )$ is equal to $|\Inv (\sigma )|$ plus a
constant that only depends on $\Des (\sigma )$; and (ii) the fact that
$\lambda _{i}$ is the number of boxes in row $i$ of the diagram.  Both
properties are preserved by rearranging the columns of the diagram and
redefining the arm-lengths $a(u)$.

Finally, in \cite[\S 5.2]{HaHaLo05}, we used (i) the fact that $\inv
(\sigma )$ enumerates certain ``inversion triples,'' plus inversions
of $\sigma $ between certain fixed pairs of boxes in row $1$ of the
diagram; and (ii) that the inversion property of a triple with respect
to the ordering $<_{2}$ in \cite[eq.~(30)]{HaHaLo05} is invariant
under a sign change $a\mapsto \overline{a}$ applied to the entry in
any one box $v$, provided that (a) for every box $u<v$ in the reading
order, $x = \sigma (u)$ satisfies $|x|>|a|$, and (b) for every box $u$
at most one row below $v$, $x=\sigma (u)$ satisfies $|x|\geq |a|$.  In
the present context, we can define inversion triples for unaugmented,
signed, possibly attacking fillings by the same rule as in \S
\ref{subsec:inv-triples}, with the modification that $\chi
_{uv}(\sigma ) = 1$ if $\sigma (u) = \sigma (v) = \overline{a}$, where
$\overline{a}$ is a negative letter.  Then property (ii) holds by
essentially the same case checking as in the original proof.
Lemma~\ref{lem:attack-pairs} holds for triples in the unaugmented
diagram $\skydg (\mu )$ with ``row $0$'' changed to ``row $1$.''  The
counterpart of Proposition~\ref{prop:inv} (the proof is the same) is
then that $\inv (\sigma )$ is the number of inversion triples plus the
number of inversions between boxes $(i,1)$, $(j,1)$, where $i<j$ and
$\mu _{i}\leq \mu _{j}$.  Hence property (i) also holds.
\end{proof}

The analog of Theorem~\ref{thm:general-H-formula} for the integral
forms $J_{\mu }(x;q,t)$, generalizing
\cite[Proposition~8.1]{HaHaLo05}, is as follows.

\begin{cor}\label{cor:general-J-formula}
Let $n(\lambda ) = \sum _{i} (i-1)\lambda _{i} $.  The integral form
Macdonald polynomials are given by the sum over unaugmented,
non-attacking fillings
\begin{equation}\label{e:J-lambda}
\begin{aligned}
J_{\lambda }&(x;q,t) = t^{n(\lambda )}D_{\mu }[X(1-t^{-1});q,t^{-1}]\\
	& = \sum _{\substack{\sigma \colon \skydg (\mu )\rightarrow [m]\\
                                  \text{non-attacking}\\}}
         x^{\sigma }q^{\maj (\sigma )} t^{n(\lambda )-\inv (\sigma )}
         \prod _{\sigma (u)=\sigma (d(u))} \left(1-q^{l(u)+1}t^{a(u)+1} \right)
         \prod _{\sigma (u)\not =\sigma (d(u))} (1-t),
\end{aligned}
\end{equation}
where $\mu $ is any rearrangement of the parts of the partition
$\lambda $.  Any boxes $u$ in row $1$ of $\skydg (\mu )$ are included
in the last factor on the second line.
\end{cor}

\begin{proof}
The first identity is Theorem~\ref{thm:general-H-formula}, combined
with the definition of $\widetilde{H}_{\lambda }$.  The proof of the
second identity is the same as the proof of
\cite[Proposition~8.1]{HaHaLo05} (see also Lemmas
\ref{lem:LLT-non-attacking} and \ref{lem:maj-inv-general}).
\end{proof}

\subsection{Stable limits of non-symmetric Macdonald polynomials}
\label{subsec:stable}

Let $(\nu ;\mu )$ denote the concatenation of two compositions, and
let $0^{m}\in \NN ^{m}$ denote the zero composition.  The diagram
$\skydg ((0^{m};\mu ))$ is just the diagram of $\mu $, shifted $m$
columns to the right.  We will identify fillings of these two diagrams
in the obvious way.  Note that $a(u) = a(u')$ for any $u\in \skydg
(\mu )$ and its corresponding box $u'\in \skydg ((0^{m};\mu ))$.

A general fact about non-symmetric Macdonald polynomials, which
follows easily from the definition, is that if $s_{i}\mu = \mu $, then
$E_{\mu }$ is $s_{i}$-invariant.  In particular, $E_{(0^{m};\mu )}$ is
symmetric in the variables $x_{1},\ldots,x_{m}$.  Setting the
remaining variables to zero therefore gives a symmetric function,
which we will denote $E_{(0^{m};\mu )}(x_{1},\ldots,x_{m};q,t)$, by
slight abuse of notation.  It is easy to see from
Theorem~\ref{thm:main} that as $m\rightarrow \infty $, this converges
to a well-defined symmetric function in infinitely many variables, the
{\it stable limit} of the non-symmetric Macdonald polynomial $E_{\mu
}$.

\begin{thm}\label{thm:stable-limit}
The stable limits of the integral form non-symmetric Macdonald
polynomials are the integral form symmetric Macdonald polynomials.
Precisely,
\begin{equation}\label{e:stable-limit}
J_{\lambda }(x_{1},\ldots,x_{m};q,t) = \Ecal _{(0^{m};\mu
)}(x_{1},\ldots,x_{m};q,t),
\end{equation}
for any $\mu $ which is a rearrangement of the parts of $\lambda $.
\end{thm}

\begin{proof}
Every non-attacking filling $\sigma $ of the unaugmented diagram
$\skydg (\mu )$ with values in $[m]$ extends to a non-attacking
augmented filling $\widehat{\sigma }$ of $\augdg ((0^{m};\mu ))$,
because the boxes $(j,0)$ for $j\leq m$ only attack other boxes in row
$0$, and the boxes $(j,0)$ for $j>m$ contain entries $\widehat{\sigma
}((j,0))>m$.  It is clear that $\maj (\widehat{\sigma }) = \maj
(\sigma )$, and that $\widehat{\sigma }$ has no descents in row $1$.
We also have $\coinv (\widehat{\sigma }) = C-\inv ({\sigma })$, where
$C = \sum _{u\in \skydg (\mu )}a(u) -|\{i<j:\mu _{i}>\mu _{j} \}|$.
Comparing \eqref{e:main-integral} and \eqref{e:J-lambda}, we see that
the result follows if $C = n(\lambda )$.  If $\mu _{1}\leq \cdots \leq
\mu _{n}$, so $\skydg (\mu )$ is the reversed diagram of the partition
$\lambda '$, this is immediate, since $n(\lambda )$ is the sum of the
arms of the boxes in the diagram of $\lambda '$.  If $\mu _{i}<\mu
_{i+1}$, then transposing these two parts of $\mu $ increases
$|\{i<j:\mu _{i}>\mu _{j} \}|$ by $1$.  On the other hand, for every
box $u\in \skydg (\mu )$, we have $a(s_{i}(u)) = a(u)$, with one
exception: for $u=(i+1,\mu _{i}+1)$, we have $a(s_{i}(u)) = a(u)+1$.
So the transposition increases $\sum _{u\in \skydg (\mu )}a(u)$ by
$1$, and $C$ remains invariant.  It follows that $C = n(\lambda )$ for
every rearrangement $\mu $ of $\lambda $.
\end{proof}

\begin{cor}\label{cor:stable-limit}
Let $\lambda ^{\circ }$ be the rearrangement of the parts of $\lambda
$ in weakly increasing order.  For any rearrangement $\mu $ of
$\lambda $, we have
\begin{equation}\label{e:monic-limit}
P _{\lambda }(x_{1},\ldots,x_{m};q,t) = \left(\frac{\prod _{u\in
\skydg (\mu )}1-q^{l(u)+1}t^{a(u)+1}}{\prod _{u\in \skydg (\lambda
^{\circ })}1-q^{l(u)}t^{a(u)+1}} \right)E _{(0^{m};\mu
)}(x_{1},\ldots,x_{m};q,t).
\end{equation}

\end{cor}

\begin{remarks}\label{rem:Knop}
(1) A related stable limit, but not truncated to the symmetric part,
is considered by Knop in a most interesting preprint \cite{Knop04},
where he conjectures an extension of the positivity theorem for
Macdonald polynomials to the non-symmetric case.

(2) The Jack polynomial version of Theorem~\ref{thm:stable-limit} is
due to Knop and Sahi \cite[Theorem~4.10]{KnSa97}.  It is possible to
prove Theorem~\ref{thm:stable-limit} directly from the orthogonality
of Macdonald polynomials, as Knop and Sahi did in the Jack case.  Then
one can reverse the arguments and deduce
Corollary~\ref{cor:general-J-formula} and
Theorem~\ref{thm:general-H-formula} from
Theorem~\ref{thm:stable-limit} and Corollary~\ref{cor:integral-form}.
In principle, this gives an alternative proof of
\cite[Theorem~2.2]{HaHaLo05}, although the only part of the original
proof that it avoids is \cite[Lemma~5.2]{HaHaLo05}.
\end{remarks}

\subsection{Symmetrization of non-symmetric Macdonald polynomials}
\label{subsec:symmetrize}

For any root system, the Macdonald polynomial $P_{\lambda }$ is the
unique monic symmetric linear combination of the polynomials $E_{\mu
}$ for $\mu $ in the Weyl group orbit of $\lambda $.  The coefficient
of $E_{\mu }$ in $P_{\mu }$ can be computed explicitly using
intertwiners \cite{Chered95}.  In the $GL_{n}$ case, this is most
conveniently expressed in terms of the polynomials $E_{\mu
}(x;q^{-1},t^{-1})$, using the fact that $P_{\mu }(x;q,t) = P_{\mu
}(x;q^{-1},t^{-1})$, as follows.

\begin{prop}[{\cite{Mac96}, \cite[Lemma 2.5(a)]{Mar99}}]\label{prop:P-from-E}
With $\lambda ^{\circ }$ as in Corollary~\ref{cor:stable-limit}, we have
\begin{equation}\label{e:P-from-E}
P_{\lambda }(x;q,t) = \prod _{u\in \skydg (\lambda ^{\circ })}
(1-q^{l(u)+1}t^{a(u)})\cdot \sum _{\mu \sim \lambda }\frac{E_{\mu
}(x;q^{-1},t^{-1})}{\prod _{u\in \skydg (\mu )}(1-q^{l(u)+1}t^{a(u)})},
\end{equation}
where the sum is over all rearrangements $\mu $ of $\lambda $.
\end{prop}

Combined with Theorem~\ref{thm:main} and
Corollary~\ref{cor:main-inverted}, we can regard \eqref{e:monic-limit}
and \eqref{e:P-from-E} as different kinds of combinatorial formulas
for $P_{\lambda }$.  In particular, letting $q,t\rightarrow 0$ in
\eqref{e:P-from-E} yields a curious combinatorial expansion for the
Schur function $s_{\lambda }(x)$, namely,
\begin{equation}\label{e:mason}
\begin{gathered}
s_{\lambda }(x) = \sum _{\mu \sim \lambda } E_{\mu }(x;\infty ,\infty
), \quad \text{where}\\
E_{\mu }(x;\infty ,\infty ) = \sum _{\substack{
\sigma \colon \skydg (\mu)\rightarrow [n]\\
\text{non-attacking}\\
\maj '(\widehat{\sigma }) = \coinv '(\widehat{\sigma })=0}} x^{\sigma }
\end{gathered}
\end{equation}
It is natural to regard the polynomials $E_{\mu }(x;\infty ,\infty )$,
which coincide with what combinatorialists refer to as {\it key
polynomials} \cite{ReinShim95}, as ``non-symmetric Schur functions.''
They are closely related but not equal to the polynomials $E_{\mu
}(x;0,0)$, which are Demazure characters \cite{Ion04}, \cite{Ion05}.
A further study of $E_{\mu }(x;\infty, \infty )$ has been made by
Sarah Mason \cite{Mason06}, who among other things found a direct,
bijective proof of \eqref{e:mason}.

\vfill

\section*{Appendix: table of $E_{\mu }(x;q,t)$}
\label{sec:tables}

The case $n=1$ is trivial.  A general closed formula for $n=2$ can be found
in Macdonald's book \cite[(6.2.7--8)]{Mac96}.  We give a short table for $n=3$.
\begin{align*}
E_{(0,0,0)} & = 1\\
E_{(1,0,0)} & = x_{1}\\
E_{(0,1,0)} & = x_{2} + {\frac{1-t}{1-q\,{t^2}}}\,x_{1}\\
E_{(0,0,1)} & = x_{3} + {\frac{1-t}{1-q\,t}}(x_{1} + x_{2})\\
E_{(1,1,0)} & = x_{1}\,x_{2}\\
E_{(1,0,1)} & = x_{1}\,x_{3} + {\frac{1-t}{1-q\,{t^2}}}\,x_{1}\,x_{2}\\
E_{(0,1,1)} & = x_{2}\,x_{3} +
                  {\frac{1-t}{1-q\,t}}(x_{1}\,x_{2}+x_{1}\,x_{3})\\  
E_{(2,0,0)} & = {{x_{1}}^2} + {\frac{q\,\left(1-t
                  \right)}{1-q\,t}}(x_{1}\,x_{2} + x_{1}\,x_{3})\\ 
E_{(0,2,0)} & = {{x_{2}}^2} + {\frac{1-t}
                                   {\left(1-q^{2}\,t^{2} \right)}}\,{{x_{1}}^2}
                  + {\frac{q\,\left( 1-t \right)}{1-q\,t}}\,x_{2}\,x_{3}\\
            & \qquad
                  + {\frac{q\,{{\left(1-t \right) }^2}}
                      {{{\left(1-q\,t \right) }}\,
                          \left(1-q^{2}\,t^{2} \right)}}\,x_{1}\,x_{3}
                  + {\frac{\left(1-t \right)\,
                            \left(1+q-q\,t-{q^2}\,{t^2} \right)}
                          {{{\left( 1-q\,t \right) }}\,
                            \left(1-q^{2}\,t^{2} \right) }}\,x_{1}\,x_{2}
\end{align*}




\providecommand{\bysame}{\leavevmode\hbox to3em{\hrulefill}\thinspace}

\end{document}